\documentclass[twoside,a4paper,11pt]{article}
\usepackage{amsmath}
\usepackage{amssymb}
\usepackage{latexsym,a4wide}
\usepackage{xy}

\input xypic
\input xy
\xyoption{all}

\newtheorem{example}{Example}[section]

\newtheorem{prop}[example]{Proposition}
\newtheorem{thm}[example]{Theorem}
\newtheorem{cor}[example]{Corollary}
\newtheorem{lem}[example]{Lemma}

\newenvironment{pf}{\noindent \textbf{Proof:}}{\rule{0em}{1ex}\hfill$\Box$\mbox{}}
\def\leq{\leqslant}
\def\geq{\geqslant}
\date{}
\begin{document}

\author{Z. Arvasi and E. Ulualan}
\title{{Homotopical Aspects of Commutative Algebras I:\\
Freeness Conditions for Crossed Squares}}
 \maketitle

\begin{abstract}
We give an alternative description of the top algebra of the free crossed
square of algebras on 2-construction data in terms of tensors and coproducts
of crossed modules of commutative algebras.
%\textbf{AMS Classification:} 13D25, 18G30, 18G55
\end{abstract}

\section*{Introduction}
Ellis \cite{ellis3} gave a description of the free crossed square of
groups of a CW-complex using topological methods. In contexts other
than groups, he has defined crossed squares, 2-crossed modules etc.
in \cite{ellis3}. For commutative algebras, 2-crossed modules have
been defined by Grandjean and Vale \cite{gv}.

Combining earlier work \cite{porter1} of Porter with  Arvasi and Porter's joint papers
\cite{ap1,ap2,ap3}, one starts to see how a study
of the links between simplicial commutative algebras and classical
constructions of homological algebra can be strengthened by
interposing crossed algebraic models for the homotopy types of
simplicial algebras. In this note, we continue this process using
these methods to give a description in terms of tensor products,
of the top corner of a free crossed square of commutative algebras
and the top term of the corresponding free 2-crossed module. The
methods are slightly different, but these results are
2-dimensional analogues of the description of the free crossed
module on a `presentation' of an algebra in terms of Kozsul-like
terms, given in \cite{porter1}.

We end with a section looking at possible links of this with
Andr\'e-Quillen homology and a squared complex form of the
contangent complex.

The results and general methods we use are inspired by those given for the
corresponding case of groups in \cite{mp}. Some of the methods of that paper,
of course, go across almost verbatim to this commutative algebra case,
but, as the audience for this paper is probably more or less disjoint from
that for \cite{mp} it seems advisable to repeat arguments from that paper in
this algebra case even when they might be safely `left to the reader'. Of
course, there are all times when the translation between the two contexts
is less easy.

\textbf{Acknowledgements}. The authors wishes to thank Professor
Timothy Porter
for his helpful comments.

\section{Preliminaries}

All algebras will be commutative and will be over the same fixed but
unspecified ground ring.

\subsection{Simplicial algebras}

Denoting the usual category of finite ordinals by $\Delta $, we obtain for
each $k\geq 0$ a subcategory $\Delta _{\leq k}$ determined by the objects $%
[j]$ of $\Delta $ with $j\leq k$. A simplicial algebra is a
functor from the opposite category $\Delta ^{op}$ to \textbf{Alg};
a $k$-truncated simplicial algebra is a functor from $(\Delta
_{\leq k}^{op})$ to \textbf{Alg}. We denote the category of
simplicial algebras by \textbf{SimpAlg} and the category of
\textit{k}-truncated simplicial algebras by
$\mathbf{Tr_{k}SimpAlg}$. By a \textit{k-truncation of a
simplicial algebra}, we mean a $k$-truncated simplicial algebra
\textbf{tr}$_{k}\mathbf{E}$ obtained by forgetting dimensions of
order $>k$ in a simplicial algebra $\mathbf{E},$ that is
restricting \textbf{E} to $\Delta _{\leq k}^{op}$. This gives a
truncation functor
\begin{equation*}
\mathbf{tr_{k}:SimpAlg\longrightarrow Tr_{k}SimpAlg}
\end{equation*}
which admits a right adjoint
\begin{equation*}
\mathbf{cosk_{k}:Tr_{k}SimpAlg\longrightarrow SimpAlg}
\end{equation*}
called the \textit{k-coskeleton functor}, and a left adjoint
\begin{equation*}
\mathbf{sk_{k}:Tr_{k}SimpAlg\longrightarrow SimpAlg},
\end{equation*}
called the \textit{k-skeleton functor}. For explicit construction
of these see \cite{dusk}. We will say that a simplicial algebra
$\mathbf{E}$ is \textit{\ k-skeletal} if the natural morphism
\textbf{sk$_{k}\mathbf{E}$}$  \rightarrow \mathbf{E}$ is an
isomorphism.

Recall that given a simplicial algebra \textbf{E}, \emph{the Moore
complex} $(\mathbf{NE},\partial )$ \emph{of} \textbf{E} is the
chain complex defined by
\begin{equation*}
(NE)_n=\bigcap_{i=0}^{n-1}\mbox{\rm Ker}d_i^n
\end{equation*}
with $\partial _n:NE_n\rightarrow NE_{n-1}$ induced from $d_n^n$ by
restriction.

The \emph{n$^{th}$ homotopy module} $\pi _{n}$(\textbf{E}) of \textbf{E} is
the n$^{th}$ homology of the Moore complex of \textbf{E}, i.e.,
\begin{equation*}
\begin{array}{rcl}
\pi _{n}(\mathbf{E}) & \cong & H_{n}(\mathbf{NE},\partial ) \\
& = & \bigcap\limits_{i=0}^{n}\mbox{\rm
Ker}d_{i}^{n}/d_{n+1}^{n+1}(\bigcap \limits_{i=0}^{n}\mbox{\rm
Ker}d_{i}^{n+1}).
\end{array}
\end{equation*}
We say that the Moore complex \textbf{NE} of a simplicial algebra is
of \emph{length} $k$ if $NE_{n}=0$ for all $n\geq k+1$ so that a
Moore complex is of length $k$ also of length $r$ for $r\geq k.$ For
example, if \textbf{E} has Moore complex of length 1, then
$(NE_{1},NE_{0},\partial _{1})$ is a crossed module and conversely.
If \textbf{NG} is of length 2, the corresponding Moore complex gives
a 2-crossed module (cf. \cite{ap3}).

\subsection{Free Simplicial Algebras}

Recall from \cite{ap1} the definition of free simplicial algebra
given by the `step-by-step' construction of Andr\'e \cite{andre}.

Let \textbf{E} be a simplicial algebra and $k\geq 1$, $k$-skeletal
be  fixed. A simplicial algebra \textbf{F} is called a
\textit{free} if

i) $F_{n}=E_{n}$ for $n<k,$

ii) $F_{k}=$ a free $E_{k}$-algebra over a set of non-degenerate
indeterminates, all of whose faces are zero except the $k^{th}$,

iii) $F_{n}$ is a free $E_{n}$-algebra over the degenerate
elements for $n>k$.

\textbf{Remark}: if \textbf{A} is a simplicial algebra, then there
exists a free simplicial algebra \textbf{E} and an epimorphism
$\mathbf{E\rightarrow A}$ which induces isomorphisms on all
homotopy modules. The details are omitted as they are
'well-known'.

\subsection{Crossed Modules of Algebras}
Throughout this paper we denote an action of $r\in R$ on $c\in C$
by $r\cdot c$.

A \emph{crossed module} is an algebra morphism $\partial
:C\rightarrow R$ with an action of $R$ on $C$ satisfying (i)
$\partial (r\cdot c)=r\partial c$ and (ii) $\partial (c)\cdot
c^{\prime }=cc^{\prime }$ for all $c,c^{\prime }\in M,r\in R$. For
the weaker notion in which condition (ii) is not required, the
models are called \textit{pre-crossed modules}.

Examples of crossed modules are:

i) Any ideal, $I$, in $R$ gives an inclusion map $I\longrightarrow
R,$ which is a crossed module. Conversely, given any crossed
$R$-module $\mu :C\longrightarrow R$, the image $I=\mu (C)$ of $C$
is an ideal in $R$.

ii) Any $R$-module $M$ can be considered as an $R$-algebra with
zero multiplication and hence $M\overset{0}{\longrightarrow }R$ is
a crossed $R$-module. Conversely, if $\mu :C\longrightarrow R$ is
a crossed $R$-module, Ker$\mu $ is an $R/\mu (C)$-module.

\section{Crossed Squares and Simplicial Algebras}
Although we will be mainly concerned with crossed squares in this
paper, some of the arguments either clearly apply or would seem to
apply in the more general case of crossed $n$-cubes and $n$-cube
complexes.

Crossed $n$-cubes in algebraic settings such as commutative
algebras, Jordan algebras, Lie algebras have been defined by
Ellis, \cite{ellis2}.

A \emph{crossed n-cube of commutative algebras} is a family of
commutative algebras, $M_A$ for $A\subseteq \langle
n\rangle=\{1,...,n\}$ together with homomorphisms $\mu
_i:M_A\rightarrow M_{A-\{i\}}$ for $i\in \langle n \rangle$ and for
$A,B\subseteq \langle n \rangle$, functions
\begin{equation*}
h:M_A\times M_B\longrightarrow M_{A\cup B}
\end{equation*}
such that for all $k\in \mathbf{k},\ a,a^{\prime }\in M_A,\
b,b^{\prime }\in M_B,\ c\in M_C,$ $i,j\in \langle n \rangle$ and
$A\subseteq B$
\begin{equation*}
\begin{array}{lrcl}
1) & \mu _ia & = & a\ \quad \mbox{\rm {\rm if}}\ i\not \in A \\
2) & \mu _i\mu _ja & = & \mu _j\mu _ia \\
3) & \mu _ih(a,b) & = & h(\mu _ia,\mu _ib) \\
4) & h(a,b) & = & h(\mu _ia,b)=h(a,\mu _ib)\ \quad \mbox{\rm {\rm if}}\ i\in
A\cap B \\
5) & h(a,a^{\prime }) & = & aa^{\prime } \\
6) & h(a,b) & = & h(b,a) \\
7) & h(a+a^{\prime },b) & = & h(a,b)+h(a^{\prime },b) \\
8) & h(a,b+b^{\prime }) & = & h(a,b)+h(a,b^{\prime }) \\
9) & k\cdot h(a,b) & = & h(k\cdot a,b)=h(a,k\cdot b) \\
10) & h(h(a,b),c) & = & h(a,h(b,c))=h(b,h(b,c)). \\
&  &  &
\end{array}%
\end{equation*}

A \emph{morphism of crossed n-cubes} is defined in the obvious way:
It is a family of commutative algebra homomorphisms, for $A\subseteq
\langle n \rangle$, $\ f_A:M_A\longrightarrow M_A^{\prime } $
commuting with the $\mu _i^{\prime}$s and $h^{\prime}s$. We thus
obtain a category of crossed $n$-cubes denoted by $\mathbf{Crs^n} $.

\vspace{0.5cm}

{\noindent \textbf{Examples.}} ${\noindent(1)}$ For $n=1,$ a crossed 1-cube
is the same as a crossed module.

{\noindent } For $n=2$ one has a crossed square:
\begin{equation*}
\xymatrix{M_{\langle 2 \rangle}\ar[r]^{\mu_2}\ar[d]_{\mu_1}&
M_{\{1\}}\ar[d]^{\mu_1}\\
M_{\{2\}}\ar[r]_{\mu_2}&M_{\emptyset}. }
\end{equation*}
Each $\mu _i$ is a crossed module as is $\mu _1\mu _2$. The
$h$-functions give actions and a pairing
\begin{equation*}
h:M_{\{1\}}\times M_{\{2\}}\longrightarrow M_{\langle 2\rangle}.
\end{equation*}
The maps $\mu _2$ (or $\mu _1)$ also define a map of crossed
modules from $(M_{\langle 2 \rangle},M_{\{2\}},\mu _1)$ to
$(M_{\langle 2 \rangle},M_\emptyset ,\mu _1)$. In fact a crossed
square can be thought of as a crossed module in the category of
crossed modules.

(2) Let $I_{1}$ and $I_{2}$ be ideals of an algebra $E$. The
commutative square diagram of inclusions;
\begin{equation*}
\xymatrix{I_1 \cap I_2 \ar[r]^-{inc.}\ar[d]_{inc.}&
I_2\ar[d]^{inc.}\\
I_1\ar[r]_{inc.}&E}
\end{equation*}
naturally comes together with actions of $E$ on $I_{1},I_{2}$ and
$I_{1}\cap I_{2}$ given by multiplication, and functions
\begin{equation*}
\begin{array}{lrll}
h: & I_{A}\times I_{B} & \longrightarrow & I_{A}\cap I_{B}=I_{A\cup B} \\
& (a,b) & \longmapsto & ab.
\end{array}%
\end{equation*}%
That this is a crossed square is easily checked.

(3) Let $\mathbf{E}$ be a simplicial algebra. Let $M(\mathbf{E},2)$
denote the following diagram
\begin{equation*}
\xymatrix{NE_2/\partial_3(NE_3)\ar[r]^-{\partial_2}\ar[d]_{\partial'_2}&NE_1\ar[d]^{\mu}\\
\overline{NE_1}\ar[r]_{\mu'}&E_1.}
\end{equation*}
Then this is the underlying square of a crossed square. The extra
structure is given as follows: $NE_1=$\text{Ker}$d_0^1$ and
$\overline{NE}_1=$\textrm{Ker}$d_1^1$. Since $E_1$ acts on
$NE_2/\partial _3NE_3,\ \overline{NE}_1$ and $NE_1,$ there are
actions of $\overline{NE}_1$ on $NE_2/\partial _3NE_3$ and $NE_1$
via $\mu ^{\prime },$ and $NE_1$ acts on $NE_2/\partial _3NE_3$ and
$\overline{NE}_1$ via $\mu .$ As $\mu $ and $\mu ^{\prime }$ are
inclusions, all actions can be given by multiplication. The $h$-map
is
\begin{equation*}
\begin{array}{rcl}
NE_1\times \overline{NE}_1 & \longrightarrow & NE_2/\partial _3NE_3 \\
(x,\overline{y}) & \longmapsto & h(x,\overline{y})=s_1x(s_1y-s_0y)+\partial
_3NE_3,%
\end{array}%
\end{equation*}
which is bilinear. Here $x$ and $y$ are in $NE_1$ as there is a
natural bijection between $NE_1$ and $\overline{NE}_1$ (by
\cite[Lemma 2.1]{ap2}). The element $\bar y$ is the image of $y$
under this. This example effectively introduces the functor
$$
\mathbf{M}(-,2):\mathbf{SimpAlg}\rightarrow \mathbf{Crs^2}.
$$
This is the case $n=2$ of a general construction of a crossed
$n$-cube from a simplicial algebra given by the first author in
\cite{a1} where the reader may find the verification of the
axioms. (This notational convention will be revisited at the end
of section \ref{sec5})

Note if we consider the above crossed square as a vertical
morphism of crossed modules, we can take its kernel and cokernel
within the category of crossed modules. In the above the morphisms
in the top left hand corner are induced from $d_2$ so
\begin{equation*}
\text{Ker}\left( \partial'_2:\frac{NE_2}{\partial _3NE_3}\longrightarrow \text{Ker%
}d_1\right) =\frac{NE_2\cap \text{Ker}d_2}{\partial _3NE_3}\cong \pi _2(%
\mathbf{E})
\end{equation*}
whilst the other map labeled $\mu $ is an inclusion so has trivial
kernel, hence the kernel of this morphism of crossed modules is
\begin{equation*}
\pi _2(\mathbf{E})\longrightarrow 0.
\end{equation*}
The image of $\partial_2$ (and $\mu' $) is an ideal in both the
algebras on the bottom line and as Ker$d_0=NE_1$ with the
corresponding Im$\mu$ being $d_2NE_2,$ the cokernel is
$NE_1/\partial _2NE_2,$ whilst $E_1/$Ker$d_0\cong E_0,$ i.e, the
cokernel of $\mu $ is $M(\mathbf{E},1)$.

In fact of course $\mu$ is not only a morphism of crossed modules,
it is a crossed module. This means that $\pi
_2(\mathbf{E})\rightarrow 0$ is in some sense a
$M(\mathbf{E},1)$-module, (cf. \cite{ae}), and that
$M(\mathbf{E},2)$ can be thought of as a crossed extension of
$M(\mathbf{E},1)$ by $\pi _2(\mathbf{E})$.

\section{ Free Crossed Squares}

Ellis, \cite{ellis3}, in 1993 presented the notion of a free
crossed square for the case of groups. In this section, we
introduce a commutative algebra version of this definition and
give a construction of a free crossed square by using the second
order Peiffer elements and the 2-skeleton of a `step-by-step'
construction of a free simplicial algebra.

We firstly adapt Ellis's definition of the free crossed square on a pair of
functions $(f_2,f_3)$ to the algebra context:

Let $\mathbf{S_{1},S_{2}}$ and $\mathbf{S_{3}}$ be sets which for
simplicity we assume are finite. Suppose given a function
$f_{2}:\mathbf{S_{2}\rightarrow }R$ from a set $\mathbf{S_{2}}$ to a
free algebra $R$ on $\mathbf{S_{1}}$. Let $\partial :M\rightarrow R$
be the free pre-crossed module on $f_{2}$. Using the action of $R$
on $M$ we can form the semi-direct product $M\rtimes R$. The
inclusion $\mu :M\rightarrow M\rtimes R$ given by $m\mapsto (m,0)$
enables us to take $M$ as an ideal of $M\rtimes R$. (Recall from
examples of crossed modules that any ideal inclusion is a crossed
module with action by multiplication.) There is also another ideal
of $M\rtimes R$ coming from $M$, namely
\begin{equation*}
\overline{M}=\{(m,r)\in M\rtimes R:\partial m=-r\}
\end{equation*}
with inclusion denoted $\bar{\mu}:\overline{M}\rightarrow M\rtimes R$.

Assume given a function from a set $\mathbf{S_{3}}$ to $M$,
$f_{3}:\mathbf{S_{3}}\rightarrow M$, which is to satisfy $\partial
f_{3}=0$. Then there is a corresponding function
$\bar{f}_{3}:\mathbf{S_{3}}\rightarrow \overline{M}$ given by
$y\mapsto (f_{3}(y),0)$.

We say a crossed square
\begin{equation*}
\xymatrix{L\ar[r]^{\partial'_2}\ar[d]_{\partial_2}&\overline{M}\ar[d]^{\mu'}\\
M\ar[r]_-{\mu}&M\rtimes R}
\end{equation*}
is \emph{totaly free} on the pair of functions $(f_2,f_3)$ if

{\noindent (i) }$(M,R,\partial )$ is the free pre-crossed module on $f_2;$

{\noindent (ii)} $\mathbf{S_3}$ is a subset of $L$ with $f_3$ and $\bar f_3$
the restrictions of $\partial _2$ and $\partial _2^{\prime }$ respectively;

{\noindent (iii)} for any crossed square
\begin{equation*}
\xymatrix{L'\ar[r]^{\tau'}\ar[d]_{\tau}&\overline{M}\ar[d]^{\mu'}\\
M\ar[r]_{\mu}&M\rtimes R}
\end{equation*}
and any function $\nu :\mathbf{S_3}\rightarrow L'$ satisfying $\tau \nu =f_3$%
, there is a unique morphism $\Phi =(\phi $,id,id,id) of crossed
squares:
$$
\xymatrix@!0{
  & L \ar[rrrr] \ar[dddd]\ar[ddl]_{\phi}
      &  & & &\overline{M} \ar[dddd]\ar@{=}[ddl]        \\
      \\
  L' \ar[rrrr]\ar[dddd]
      &  &  & &\overline{M} \ar[dddd] \\
      \\
  & M \ar[rrrr]\ar@{=}[ddl]
      &  & &   &M\rtimes R \ar@{=}[ddl]               \\
      \\
  M \ar[rrrr]
      &  &  &  & M\rtimes R    }
$$
such that $\phi \nu ^{\prime }=\nu ,$ where $\nu ^{\prime
}:\mathbf{S_3}\to L $ is the inclusion.

We denote such a free crossed square of algebras by
$(L,M,\overline{M},M\rtimes R). $

We know the free pre-crossed module on
$f:\mathbf{S_{2}}\rightarrow R$ is $\partial
:R^{+}[\mathbf{S_{2}}]\rightarrow R$, so the function
$$f_{3}:\mathbf{S_{3}}\rightarrow M,\ \ (\text{with}\ \ M=R^{+}[\mathbf{S_{2}}])$$
is precisely the data $(\mathbf{S_{3}},f_{3})$ for 2-dimensional
construction data in the simplicial context. We thus need to recall
the 2-dimensional construction data in a free simplicial algebra
(cf. \cite{ap2}). This 2-dimensional form can be pictured by the
diagram
\begin{equation*}
\xymatrix@C=3pc{\mathbf{E^{(2)}}:\cdots (R[s_0 (\mathbf{S_2}),s_1
\mathbf{S_2}])[\mathbf{S_3}]\ar@<2ex>[r]^-{d_0,d_1,d_2}\ar@<1ex>[r]\ar@<0pt>[r]
&
R[\mathbf{S_2}]\ar@<1ex>[r]^-{d_0,d_1}\ar@<0ex>[r]\ar@<1ex>[l]\ar@<2ex>[l]^-{s_0,s_1}&R
\ar@<1ex>[l]^-{s_0}}
\end{equation*}
with the simplicial morphisms given as in \cite{ap3}. Here
$\mathbf{S_{2}} =\{S_{1},\ldots ,S_{n}\}$ and
$\mathbf{S_{3}}=\{S_{1}^{\prime },\ldots ,S_{m}^{\prime }\}$ are
finite sets with $\pi _{0}(\mathbf{E}^{(2)})\cong
B=R/(t_{1},\ldots ,t_{n})$ as an $R$-algebra where $t_i=\partial
S_i$.

\subsection{Free crossed squares exist}

\begin{thm}
A totally free crossed square $(L,M,\overline{M},M\rtimes R)$ exists
on the 2-dimensional construction data and is given by
$\mathbf{M(E^{(2)},2)}$ where $\mathbf{E^{(2)}}$ is the 2-skeletal
free simplicial algebra defined by the construction data.
\end{thm}

\begin{pf}
Suppose given the 2-dimensional construction data for a free
simplicial algebra, \textbf{E}, which we will take as above as the
data for a totally free crossed square. We will not assume
detailed knowledge of \cite{ap1} so we start with $R$ and $\
f_{2}:\mathbf{S_{2}}\rightarrow R$ and form$\
M=R^{+}[\mathbf{S_{2}}]=(S_{1},\ldots ,S_{n}).$ This gives
$\partial _{1}:R^{+}[\mathbf{S_{2}}]\rightarrow R$ as free
pre-crossed module on $f_{2}.$ The semidirect product gives us
back
\begin{equation*}
R[\mathbf{S_{2}}]\cong M\rtimes R
\end{equation*}
and we can identify this with $\mathbf{E_{1}^{(2)}}$. This identification
also makes
\begin{equation*}
M\cong \mbox{\rm Ker}d_{0}^{1}
\end{equation*}
for the $d_{0}^{1}$ of $\mathbf{E^{(2)}}$.

Next form $\overline{M}=\{(m,r)\in M\rtimes R:\partial m=-r\}$. As
$m\in R^{+}[\mathbf{S_{2}}]$, writing $m=\sum r_{\alpha }S^{\alpha
}$ for multi-indices $\alpha $, we get $\partial m=\sum r_{\alpha
}t^{\alpha }$ where $t_{i}=f_{2}(S_{i})$. Thus we can identify
$\overline{M}$ with $(S_{1}-t_{1},\ldots ,S_{n}-t_{n})$ which is
exactly Ker$d_{1}^{1}$ (for this see \cite{ap1}).

Now $\ f_{3}:\mathbf{S_{3}}\rightarrow \mbox{\rm Ker}\partial
_{1}= \mbox{\rm Ker}(\partial :NE_{1}^{(2)}\rightarrow
NE_{0}^{(2)})\subset R^{+}[\mathbf{ S_{2}}].$ We know that this
allows us to construct $\mathbf{E_{2}^{(2)}}$, and hence
$\mathbf{E_{n}^{(2)}}$ for $n\geq 3$, and in addition that taking
\begin{equation*}
L=NE_{2}^{(2)}/\partial _{3}(NE_{3}^{(2)}),
\end{equation*}
gives a crossed square
\begin{equation*}
\xymatrix{L\ar[r]^{\partial'_2}\ar[d]_{\partial_2}&\overline{M}\ar[d]^{\lambda'}\\
M\ar[r]_{\lambda}&E^{(2)}_1}
\end{equation*}
which is $M(\mathbf{E}^{(2)},2)$. We claim this is the totally free crossed
square on the construction data.

At this stage it is worth nothing that there would seem to be no simple
adjointness statement between $\mathbf{M(-,2)}$ and some functor that would
give a quick proof of freeness. The problem being that $\mathbf{M(-,2)}$
seems to be an adjoint only up to some sort of coherent homotopy. To avoid
this difficulty we use a more combinatorial approach involving the higher
order Peiffer elements and a more concrete description of $L$.

In \cite{ap1}, the first author and Porter analysed in general the
structure of algebras of boundaries such as $\partial
_{3}(NE_{3}^{(2)}).$ There they showed that $NE_{1}^{(2)}$ is
generated as an ideal by elements of the following forms:

For all $x\in NE_{1}^{(2)},~y\in NE_{2}^{(2)}$,
\begin{align*}
C_{(1,0)(2)}(x\otimes y) & =  (s_{1}s_{0}x-s_{2}s_{0}x)s_{2}y, \\
C_{(2,0)(1)}(x\otimes y) & =  (s_{2}s_{0}x-s_{2}s_{1}x)(s_{1}y-s_{2}y), \\
C_{(2,1)(0)}(x\otimes y) & = s_{2}s_{1}x(s_{0}y-s_{1}y+s_{2}y);\\
\intertext{ whilst for all $x,~y\in NE_{2}$,}
C_{(1)(0)}(x\otimes y) & =  s_{1}x(s_{0}y-s_{1}y)+s_{2}(xy), \\
C_{(2)(0)}(x\otimes y) & =  (s_{2}x)(s_{0}y), \\
C_{(2)(1)}(x\otimes y) & =  s_{2}x(s_{1}y-s_{2}y).
\end{align*}
We know that $\partial _{3}(NE_{3}^{(2)})$ is generated by elements
of the forms
\begin{equation*}
(s_{1}s_{0}d_{1}S_{i}-s_{0}S_{i})S_{j}^{\prime },\ \
(s_{0}S_{i}-s_{1}S_{i})(s_{1}d_{2}S_{j}^{\prime }-S_{j}^{\prime }),\ \
s_{1}S_{i}(s_{0}d_{2}S_{j}^{\prime }-s_{1}d_{2}S_{j}^{\prime }+S_{j}^{\prime
}),
\end{equation*}
and for $S_{i}^{\prime },$ $S_{j}^{\prime }\in \mathbf{S_{2}}$,
\begin{equation*}
S_{i}^{\prime }(s_{1}d_{2}S_{j}^{\prime }-S_{j}^{\prime }),\ \ S_{i}^{\prime
}(S_{j}^{\prime }+s_{0}d_{2}S_{j}^{\prime }-s_{1}d_{2}S_{j}^{\prime }),\ \
(s_{0}d_{2}S_{i}^{\prime }-s_{1}d_{2}S_{i}^{\prime }+S_{i}^{\prime
})(s_{1}d_{2}S_{j}^{\prime }-S_{j}^{\prime }),
\end{equation*}
which are the second order Peiffer elements defined in \cite{ap1},
where $S_{i}\in NE_{1}=\mbox{\rm Ker}d_{0}=R^{+}[\mathbf{S_{2}}]$
and $S_{i}^{\prime }\in
NE_{2}=R[s_{0}(\mathbf{S_{2}})]^{+}[s_{1}(\mathbf{S_{2}}),\mathbf{S_{3}}]\cap
(s_{0}(\mathbf{S_{2}})-s_{1}(\mathbf{S_{2}})).$

The above diagram can thus be realised as
\begin{equation*}
\xymatrix{\dfrac{R[s_0({\bf S_2})]^+[s_1({\bf S_2}),{\bf S_3}]\cap
(s_0({\bf S_2})-s_1({\bf
S_2}))}{P_2}\ar[r]^-{\partial'_2}\ar[d]_-{\partial_2}&\overline{R^+
[{\bf S_2}]}
\ar[d]^{\lambda'}\\
R^{+}[\mathbf{S_2}]\ar[r]_{\lambda}&R[\mathbf{S_2}]}
\begin{array}{cc}
&\\
&\\
&(\ast)\\
\end{array}
\end{equation*}
where $P_{2}$ is the second order Peiffer ideal which is in fact
just $\partial _{3}(NE_{3}^{(2)}).$

Given any crossed square $(T,M,\overline{M},M\rtimes R)$ with a
function $\nu :\mathbf{S_{3}}\rightarrow T$, then there exists a
morphism
\begin{equation*}
\phi :(L,M,\overline{M},M\rtimes R)\longrightarrow
(T,M,\overline{M},M\rtimes R)
\end{equation*}
given by
\begin{equation*}
\phi (S_{i}^{\prime }+P_{2})=\nu (S_{i}^{\prime })
\end{equation*}
such that $\phi \nu ^{\prime }=\nu $ where $\nu
:\mathbf{S_{3}}\rightarrow L$ is a function. The existence of
$\phi $ follows by using the freeness property of the algebra
$NE_{2}^{(2)}$ and then restricting to
$R[s_{0}(\mathbf{S_{2}})]^{+}[s_{1}(\mathbf{S_{2}}),
\mathbf{S_{3}}]\cap
(s_{0}(\mathbf{S_{2}})-s_{1}(\mathbf{S_{2}})).$ The ideal
generating elements of $P_{2}$ are then easily shown to have
trivial image in $T$ as that algebra is part of the second crossed
square.

Thus diagram $(\ast )$ is the desired totally free crossed square on
the 2-dimensional construction data. The crossed square properties
of $(L,M,\overline{M},M\rtimes R)$ may be easily verified or derived
from the fact that this is exactly $\mathbf{M(E}^{(\mathbf{2})},2).$
\end{pf}

\vspace{0.5cm}

\textbf{Remark}:

 At this stage, it is important to note that nowhere
in the argument was use made of the freeness of the 1-skeleton. If
\textbf{E} is any 1-skeletal simplicial algebra and we form a new
simplicial algebra \textbf{F} by adding in a set $\mathbf{S}_{3}$
of new generators in dimension 2, so that for instance,
$F_{2}=E_{2}^{+}[\mathbf{S}_{3}],$ a free algebra on
$\mathbf{S}_{3}$, then we can use \thinspace
$M=NE_{1}=$Ker$d_{0}^{1}$ as before even though it need not be
free. The corresponding $\overline{M}$ is then isomorphic to
Ker$d_{1}^{1}$ with the bottom right hand corner being $E_{1}$.
The `construction data' is now replaced by data for killing some
elements of $\pi _{1}(E),$ specified by
$f_{3}:\mathbf{S}_{3}\rightarrow M.$ We introduce the term
`totally free crossed square' for the type of free crossed square
constructed in the above theorem, using free crossed square for
the more general situation in which $(M,E,\partial )$ and $f_{3}$
are specified and no requirement $(M,E,\partial )$ to be a free
pre-crossed module is made.

\subsection{ The $n$-type of the $k$-skeleton}

As in the other paper in this series, we will use the
`step-by-step' construction of a free simplicial algebra to
observe the way in which the models react to the various steps of
the construction.

By a `step-by-step' construction of a free simplicial algebra,
there are simplicial inclusions
\begin{equation*}
\mathbf{E}^{(0)}\subseteq \mathbf{E}^{(1)}\subseteq
\mathbf{E}^{(2)}\subseteq \cdots
\end{equation*}%
The functor, $\mathbf{M(\ ,\ }n\mathbf{)}$, from the category of
simplicial algebras to that of crossed $n$-cubes gives the
corresponding inclusions
\begin{equation*}
\mathbf{M}(\mathbf{E}^{(0)},\ n)\hookrightarrow
\mathbf{M}(\mathbf{E}^{(1)},\ n)\hookrightarrow
\mathbf{M}(\mathbf{E}^{(2)},\ n)\hookrightarrow \mathbf{\cdots }
\end{equation*}%
We investigate $\mathbf{M(E^{(i)}},n)$, for $n=0,1,2$, and varying
$i$.

Firstly look at $\mathbf{M(E^{(0)}},n),$where the 0-skeleton
$\mathbf{E}^{(0)}\,$ is
\begin{equation*}
\begin{array}{lccc}
\mathbf{E}^{(0)}: & \cdots \longrightarrow R\longrightarrow R\longrightarrow
R & \overset{f}{\longrightarrow } & B%
\end{array}%
\end{equation*}
with the $d_i^n=s_j^n=\ $identity homomorphisms.

For $n=0,\,$there is an equality
$\mathbf{M(E^{(0)}},0)=E_0^{(0)}/d_1(\mbox{\rm Ker}d_0)=R, $
 and so $\mathbf{M(E^{(0)},}0)$ is just an algebra of $0$-simplices of \textbf{E}.

For $n=1, \mathbf{M(E^{(0)}},1)$ is $NE_1^{(0)}/\partial
_2NE_2^{(0)}\rightarrow E_0. $ It is easy to show that
$NE_1^{(0)}/\partial _2NE_2^{(0)}$ is trivial in the 0-skeleton
$\mathbf{E}^{(0)}$ and hence
\begin{equation*}
\mathbf{M(E^{(0)}},1)\cong (0\longrightarrow R).
\end{equation*}

For $n=2$, \ $\mathbf{M(E^{(0)}},2)$ is the trivial crossed square
\begin{align*}
\begin{aligned}\xymatrix{ NE_2 / d^{3}_{3} (NE_3 ) \ar[d]\ar[r]&\text{Ker}d_0
   ^1\ar[d]\\
   \text{Ker}d_1^1\ar[r]&E_1}\end{aligned}&= \begin{aligned}\xymatrix{0 \ar[d]\ar[r]&0\ar[d]\\
   0\ar[r]&R.}\end{aligned}
\end{align*}
Next take $\mathbf{M(E^{(1)}},n)$ and recall that the 1-skeleton
\textbf{E}$^{(1)}$ is
\begin{equation*}
\xymatrix@C=3pc{\mathbf{E^{(1)}}:\cdots R[s_0 (\mathbf{S_2}),s_1(
\mathbf{S_2})]\ar@<2ex>[r]^-{d_0,d_1,d_2}\ar@<1ex>[r]\ar@<0pt>[r] &
R[\mathbf{S_2}]\ar@<1ex>[r]^-{d_0,d_1}\ar@<0ex>[r]\ar@<1ex>[l]\ar@<2ex>[l]^-{s_0,s_1}&R
\ar@<1ex>[l]^-{s_0}\ar[r]^{f}&R/I.}
\end{equation*}
For $n=0$, it follows that $\mathbf{M(E^{(1)},} 0)$ is $
E_0^{(1)}/d_1(\mbox{\rm Ker}d_0)\cong R/I $ which is $\pi
_0(\mathbf{E} ^{(1)})\cong\pi _0(\mathbf{E})$.

Let $n=1.$ We have that%
\begin{align*}
\mathbf{M(E}^{(1)}\mathbf{,\ }1\mathbf{)} & = NE_{1}^{(1)}/\partial
_{2}NE_{2}^{(1)}\rightarrow E_{0}^{(1)} \\
& =  R^{+}[\mathbf{S_{2}}]/P_{1}\rightarrow R%
\end{align*}
which is the free crossed module. In fact this is the free crossed
module on the (generalised) presentation
$(\mathbf{S_{1};S_{2}},f_{2})$. As pointed out in \cite{ap1}, it is
often convenient to generalise the notion of a presentation
$P=(R\,;x_{1},\ldots ,x_{n})$ of an $R$-algebra $B$ in this way and
$\mathbf{E}^{(1)}$ is the 1-skeleton of the free simplicial algebra
generated by this presentation, then
\begin{equation*}
\delta :NE_{1}^{(1)}/\partial _{2}(NE_{2}^{(1)})\longrightarrow NE_{0}^{(1)}
\end{equation*}%
is the free crossed module on $\{S_{1},\ldots ,S_{n}\}\rightarrow
R$. This has a neat description (cf. \cite{porter1}) as follows: The
Peiffer ideal $P_1 =\partial_2(NE^{(1)}_2)$ contains all such terms
as $S_iS_j-\delta(S_i)S_j$, so any polynomial in the $S_i$'s can be
reduced mod $P_1$ to a linear form, hence each coset has a
representative of the form $\sum r_iS_i$. As $S_iS_j=S_jS_i$, these
representatives are nonunique and so the free crossed $C$ is $R^n$
factored by all $\delta(S_i)S_j-\delta(S_j)S_i$ we thus have
\begin{equation*}
\pi _{1}(\mathbf{M}(\mathbf{E}^{(1)},1))\cong \mathrm{Ker}(C\longrightarrow
R)\cong H_{2}(B,B)
\end{equation*}
the second Andr\'{e}-Quillen homology group, where $C\cong
R^{n}/\mathrm{Im}d, $ for $d: \Lambda ^{2}R^{n}\rightarrow R$, the
second Kozsul differential, see \cite{porter1} for details. Thus
\begin{align*}
\pi _{0}(\mathbf{M}(\mathbf{E}^{(1)},1)) & \cong  B \\
\pi _{1}(\mathbf{M}(\mathbf{E}^{(1)},1)) & \cong  H_{2}(B,B)
\intertext{whilst} \pi _{i}(\mathbf{M}(\mathbf{E}^{(1)},1))&\cong
0.
\end{align*}

For $n=2$,
$$
 NE_{2}^{(1)}=(R[s_{0}(\mathbf{S_{2}})]^{+}[s_{1}(\mathbf{S_{2}})])\cap
(s_{0}(\mathbf{S_{2}})-s_{1}(\mathbf{S_{2}})),
$$

$\mathbf{M(E}^{(1)},2)$ simplifies to give (up to isomorphism)
\begin{align*}
\begin{aligned}
   \xymatrix{ NE_2 / d^{3}_{3} (NE_3 ) \ar[d]\ar[r]&\text{Ker}d_0
   ^1\ar[d]\\
   \text{Ker}d_1^1\ar[r]&E_1}
   \end{aligned}&=
\begin{aligned}
   \xymatrix{\dfrac{(R[s_0({\bf S_2})]^+[s_1({\bf S_2})])\cap
(s_0({\bf S_2})-s_1({\bf S_2}))}{P_2}
\ar[d]\ar[r]&\overline{R^+[{\bf
S_2}]}\ar[d]\\
   R^+[{\bf S_2}]\ar[r]&R[{\bf S_2}]}
   \end{aligned}
  \end{align*}
which is a crossed square.

Let us next look at $\mathbf{M(E}^{(2)},n).$ Recall the 2-skeleton
$\mathbf{E}^{(2)}$
\begin{equation*}
\xymatrix@C=3pc{(R[s_0 (\mathbf{S_2}),s_1
\mathbf{S_2}])[\mathbf{S_3}]\ar@<2ex>[r]^-{d_0,d_1,d_2}\ar@<1ex>[r]\ar@<0pt>[r]
&
R[\mathbf{S_2}]\ar@<1ex>[r]^-{d_0,d_1}\ar@<0ex>[r]\ar@<1ex>[l]\ar@<2ex>[l]^-{s_0,s_1}&R
\ar@<1ex>[l]^-{s_0}\ar[r]^f&R/I.}
\end{equation*}
The following equalities can be easily obtained by direct
calculation: \thinspace for $n=0,$
\begin{equation*}
\mathbf{M(E^{(2)}},0)=E_{0}/d_{1}(\mbox{\rm Ker} d_{0})\cong \pi
_{0}(\mathbf{E}^{(2)})=\mathbf{M(E}^{(1)},0).
\end{equation*}
For $n=1,$
\begin{equation*}
\mathbf{M(E}^{(2)},1)\cong (R^{+}[\mathbf{S_{2}}
]/P_{1}\rightarrow R)=\mathbf{M(E}^{(1)},1).
\end{equation*}%
and there is an isomorphism
\begin{equation*}
\pi _{2}(\mathbf{E}^{(2)})\cong \mbox{\rm Ker}\left(
NE_{2}^{(2)}/\partial _{3}(NE_{3}^{(2)})\longrightarrow
E_{1}^{(2)}\right).
\end{equation*}
 Finally, let $n=2.$ Since by an earlier result of this section,
 $\mathbf{M(E}^{(2)},2)$ corresponds to the free
crossed square, we obtain:
\begin{align*}
\begin{aligned}
   \xymatrix{ NE_2^{(2)} / d^{3}_{3} (NE_3^{(2)} ) \ar[d]\ar[r]&\text{Ker}d_0
   ^1\ar[d]\\
   \text{Ker}d_1^1\ar[r]&E_1}
         \end{aligned}&=
         \begin{aligned}
    \xymatrix{\dfrac{R[s_0({\bf S_2})]^+[s_1({\bf S_2}),\bf{S_3}]\cap
(s_0({\bf S_2})-s_1({\bf S_2}))}{P_2}
\ar[d]\ar[r]&\overline{R^+[{\bf
S_2}]}\ar[d]\\
   R^+[{\bf S_2}]\ar[r]&R[{\bf S_2}].} \\
  \end{aligned}
  \end{align*}
This reduces to the earlier case if $\mathbf{S_3}$ is empty. Thus we
have the following relations
$$
\mathbf{M(E}^{(2)},0)=\mathbf{M(E}^{(1)},0), \quad
\mathbf{M(E}^{(2)},1)=\mathbf{M(E}^{(1)},1)
$$
but $\mathbf{M(E}^{(2)},2)$ and $\mathbf{M(E}^{(3)},2)$ need not to
be the same due to the additional influence of $\mathbf{S_3}$. Of
course it is clear that, in general:
$$
\mathbf{M(E}^{(i)},n)=\mathbf{M(E}^{(i+1)},n) \quad \text{if}\ i\geq
n+1.
$$
Clearly these top left hand corner terms are unwieldy to handle
and we will seek in section \ref{sec5} an alternative description.

\section{Squared Complexes}

The first author and M. Ko\c{c}ak defined $n$-crossed complexes of
algebras in \cite{ak} as the analogue for commutative algebras of
the notion introduced by Ellis, \cite{ellis3}, in homotopy theory.
In this paper we will only need the case $n=2$, which we shall
call a \emph{squared complex}; it consists of a diagram of algebra
homomorphisms
$$
\xymatrix{(*)&&\cdots\ar[r]&C_5\ar[r]^{\partial_5}&C_4\ar[r]^{\partial_4}&L\ar[r]^{\lambda'}\ar[d]_{\lambda}&N
\ar[d]^{\mu'}\\
&&&&&M\ar[r]_{\mu}&R}
$$
together with an action of $R$ on $L,N,M$ and $C_i$ for $i\geq 4,$
and a function $h:M\times N\rightarrow L.$ The following axioms need
to be satisfied.

$(i)$ The square
\begin{equation*}
\left(
\begin{array}{cc}
  \xymatrix{L\ar[r]\ar[d]&N\ar[d]\\
 M\ar[r]&R}
\end{array}
\right)
\end{equation*}
is a crossed square.

$(ii)$  $C_{n}$ is an $A$-module for $n\geq 4$ with $A=R/\left\{ \mu
(M)+\mu ^{\prime }(N)\right\} .$

$(iii)$ The action of $R$ on $C_{n},n\geq 4,$ is such that $\mu (M)$
and $\mu ^{\prime }(N)$ operate trivially. Thus each $C_{n}$ is an
$A$-module.

$(iv)$ each $\partial _{n}$ is $A$-module homomorphism and for
$n\geq 4$, $\partial _{n}\partial _{n+1}=0$.

A \emph{morphism} of square complexes
\begin{equation*}
\phi :\left(
\begin{array}{cc}
C_i, & \left(
\begin{array}{cc}
L & N \\
M & R%
\end{array}
\right)%
\end{array}
\right) \longrightarrow \left(
\begin{array}{cc}
C_i^{\prime }, & \left(
\begin{array}{cc}
L^{\prime } & N^{\prime } \\
M^{\prime } & R^{\prime }
\end{array}
\right)%
\end{array}
\right)
\end{equation*}
consists of a morphism of crossed squares $(\phi _L,\phi _M,\phi
_N,\phi _R)$ together with a family of $\phi _R$-equivariant
homomorphisms $\phi _i$,  $i\geq 4$ satisfying $\phi _L\partial
_4=\partial _4^{\prime }\phi _4$ and $\phi _{i-1}\partial
_i=\partial _i^{\prime }\phi _i$ for $i\geq 5.$ There is clearly a
category \textbf{SqComp} of squared complexes.

By a (totally) free squared complex, we will mean one in which the
crossed square is (totally) free, and in which each $C_n$ is free
as a $\pi _0$-module for $i\geq 3.$

\begin{prop}
There is a functor
\begin{equation*}
C(\ ,2):\mathbf{SimpAlg} \longrightarrow \mathbf{SqComp}
\end{equation*}
such that free simplicial algebras are sent to totally free
squared complexes.
\end{prop}

\begin{pf}
Let $\mathbf{E}$ be a simplicial algebra. We will define a squared
complex $ C(\mathbf{E},2)$ by specifying $C(\mathbf{E},2)_{A}$ for
each $A\subseteq \langle 2\rangle$ and for $n\geq
3,C(\mathbf{E},2)_{n}.$ As usual, (cf. the other paper in this
series {\cite{ap1, ap2, ap3}), we will denote by $D_{n}$ the ideal
of $NE_{n}$ generated by degenerate elements. }

For $A\subseteq \langle 2\rangle$, we define in particular
\begin{equation*}
C(\mathbf{E},2)_{\langle
2\rangle}=\mathbf{M}(\mathbf{sk}_{2}\mathbf{E},2)_{\langle
2\rangle}=\dfrac{NE_2}{\partial_3(NE_3 \cap D_3)}.
\end{equation*}%
We do not need to define $\mu _{i}$ and the $h$-maps relative to
these algebras as they are already defined in the crossed square
$\mathbf{M}(\mathbf{sk}_{2} \mathbf{E},2)_{A}.$

For $n\geq 3,$ we set
\begin{equation*}
C(\mathbf{E},2)_{n}=\frac{NE_{n}}{(NE_{n}\cap D_{n})+d_{n+1}(NE_{n+1}\cap
D_{n+1})}.
\end{equation*}
As this is part of the crossed complex associated to \textbf{E},
we can take the structure maps to be those of that crossed
complex, cf. \cite{ap1}. The terms are all modules over the
corresponding, $\pi _{0}$ as is easily checked. The final missing
piece, $\partial _{3},$ of the structure is induced by the
differential $\partial _{3}$ of $NE.$

The axioms for a squared complex can now be verified the known
results for crossed squares and for crossed complexes with a direct
verificatiton of those axioms relating to the interaction of the two
parts of the structure, much as in \cite{ap1}.

Now suppose the simplicial algebra is free. The proof of the
freeness of $\mathbf{M}(sk_{2}\mathbf{E},2)$ together with the
freeness of the crossed complex of a free simplicial algebra,
\cite{ap1}, now completes the proof.
\end{pf}

Suppose that $\rho $ is a general squared complex. The \emph{homotopy modules%
} $\pi _{n}(\rho ),n\geq 0$ of $\rho $ are defined in \cite{ak} to
be the homology modules of the complex
\begin{equation*}
\xymatrix{\cdots\ar[r]^{\partial_6}&C_5\ar[r]^{\partial_5}&C_4
\ar[r]^{\partial_4}& L\ar[r]^-{\partial_3}&M\rtimes
N\ar[r]^-{\partial_2}\ar[r]&0}
\end{equation*}%
with $\partial _{3}(l)=(-\lambda ^{\prime }l,\lambda l)$ and
$\partial _{2}(m,n)=\mu (m)+\mu ^{\prime }(n)$. The axioms of a
crossed square guarantee (see \cite{ap3}) that $\partial _{3}$ and
$\partial _{2}$ are homomorphisms with $\partial _{4}(C_{4})$ an ideal in Ker$(\partial _{3})$, \ $%
\partial _{3}(L)$ an ideal in Ker$(\partial _{2})$, and $\partial _{2}(M\rtimes
N)$ an ideal in $R$.

\begin{prop}
The homotopy groups of $C(\mathbf{E},2)$ are isomorphic to those of $\mathbf{%
E}$ itself.
\end{prop}

\begin{pf}
Again this is a consequence of well-known results on the two parts of the
structure.
\end{pf}

\section{Alternative Description of Freeness \label{sec5}}

In the context of CW-compexes, Ellis \cite{ellis3} gave a neat
description of the top algebra $L$ in (totally) free crossed
squares.  A free simplicial algebra is the algebraic analogue of a
CW-complex so one would expect a similar result to hold in this
setting. For this we need two constructions.

\subsection{Tensor Products}

Suppose that $\mu :M\rightarrow R$ and $\nu :N\rightarrow R$ are crossed
modules of commutative algebras over $R$. The algebras $M$ and $N$ act on
each other, and themselves, via the action of $R$. The \emph{tensor product}
$M\otimes N$ is the algebra generated by the symbols $m\otimes n$ for $m\in
M,\ n\in N$ and $r\in R$ subsect to the relations%
\begin{equation*}
\begin{array}{rrll}
\text{(i)} & r\left( m\otimes n\right) & = & rm\otimes n=m\otimes rn \\
\text{(ii)} & \left( m+m^{\prime }\right) \otimes n & = & m\otimes
n+m^{\prime }\otimes n \\
& m\otimes (n+n^{\prime }) & = & m\otimes n+m\otimes n^{\prime } \\
\text{(iii)} & \left( m\otimes n\right) \left( m^{\prime }\otimes n^{\prime
}\right) & = & \left( mm^{\prime }\otimes nn^{\prime }\right)%
\end{array}
\end{equation*}
where $m^{\prime }\in M$ and $n^{\prime }\in N.$ There are
morphisms $ \lambda :M{\otimes }N\rightarrow M,\ m{\otimes
}n\mapsto m\cdot n=m\nu (n)$ and $\lambda ^{\prime }:M{\otimes
}N\rightarrow N,\ m{\otimes }n\mapsto n\cdot m=\mu (m)n$. The
algebra $R$ acts on $M\otimes N$ by $r\cdot (m\otimes n)=r\cdot
m\otimes n=m\otimes r\cdot n,$ and there is a function $ h:M\times
N\rightarrow M\otimes N,$ \ $(m,n)\mapsto m{\otimes }n$. It is
verified in \cite{ak} that this structure is a crossed square
\begin{equation*}
\xymatrix{M\otimes N \ar[r]^-{\lambda}\ar[d]_{\lambda'}&N\ar[d]^{\nu}\\
M\ar[r]_{\mu}&R}
\end{equation*}%
with the universal property of extending the corner
\begin{equation*}
\xymatrix{&N\ar[d]^{\nu}\\
M\ar[r]_{\mu}&R.}
\end{equation*}

\subsection{Coproducts}

The following construction is due to Shammu \cite{nizar}.

Let $(M,R,\partial _{1}),(N,R,\partial _{2})$ be crossed
$R$-modules. Then $N $ acts on $M$, and $M$ acts on $N$, via the
given actions of $R$. Let $ M\rtimes N$ denote the semidirect
product with the multiplication given by
\begin{equation*}
(m,n)(m^{\prime },n^{\prime })=(mm^{\prime },n\partial _{2}(m^{\prime
})+\partial _{2}(m)n^{\prime }+nn^{\prime })
\end{equation*}
and injections
\begin{equation*}
\begin{array}{ccccc}
i^{\prime }\colon & M\rightarrow M\rtimes N & \mbox{ and } & j^{\prime
}\colon & N\rightarrow M\rtimes N \\
& m\mapsto (m,0) &  &  & n\mapsto (0,n).
\end{array}%
\end{equation*}
We define the pre-crossed module
\begin{equation*}
\begin{array}{cccl}
\underline{\delta }\colon & M\rtimes N & \longrightarrow & R \\
& (m,n) & \longmapsto & \partial _{1}(m)+\partial _{2}(n).
\end{array}
\end{equation*}
Let $P$ be the ideal of $M\rtimes N$ generated by elements of the
form
\begin{equation*}
(m,n)(m^{\prime }n^{\prime })-\underline{\delta }(m,n)(m^{\prime },n^{\prime
})=(-\partial _{1}(m)n,m\partial _{2}(n))
\end{equation*}
for all $(m,n),(m^{\prime },n^{\prime })\in M\rtimes N.$ Thus we are able to
form the quotient algebra $M\rtimes N/P$ and obtain an induced morphism
\begin{equation*}
\partial \colon M\rtimes N/P\longrightarrow R
\end{equation*}%
given by
\begin{equation*}
\partial (m,n)+P=\partial _{1}m+\partial _{2}n.
\end{equation*}%
Let $q\colon M\rtimes N\rightarrow M\rtimes N/P$ be projection and
let $ i=qi^{\prime },j=qj^{\prime }.$ Then $M\sqcup N=M\rtimes
N/P,$ with, the morphism $i,j,$ is the \emph{coproduct } of $M,N$
in the category $\mathbf{ XMod_{k}}$. The above notation can be
summarised in the following diagram:
\begin{equation*}
\xymatrix{M\ar[dr]_{i}\ar[dd]_{\partial_1}\ar[rr]^{i'}&&M\rtimes
N\ar[dl]^{q}\\
&M\rtimes N/P \ar[dl]^{\partial}&\\
R&&N.\ar[uu]^{j'}\ar[ll]^{\partial_2}\ar[ul]^{j}}
\end{equation*}

\begin{prop}
Let
\begin{equation*}
\left(
\begin{array}{cc}
\xymatrix{L\ar[r]\ar[d]&\overline{M}\ar[d]\\
M\ar[r]&M\rtimes R}
\end{array}%
\right)
\end{equation*}
be a free crossed square on the 2-dimensional construction data or
on functions $(f_{2},f_{3})$ as described above. Let $\partial
:C\rightarrow M\rtimes R$ be the free crossed module on the
function \textbf{S}$_{3}\rightarrow M\rtimes R$ given by $y\mapsto
(f_{3}y,0)$. Form the crossed module $\partial ^{\prime }:M\otimes
\overline{M}\rightarrow M\rtimes R$, then
\begin{equation*}
L\cong \{(M\otimes \overline{M})\sqcup C\}/\thicksim
\end{equation*}
where $\thicksim$ corresponds to the relations
\begin{equation*}
\begin{array}{lllll}
1) &  & i_{M\otimes \overline{M}}(\partial c\otimes \overline{n}) & \thicksim & j(c)-j(%
\overline{n}\cdot c) \\
2) &  & i_{M\otimes \overline{M}}(m\otimes \partial c) & \thicksim & j(m\cdot c)-j(c)%
\end{array}%
\end{equation*}%
for $c\in C,m\in M$ and $\overline{n}\in \overline{M.}$

The homomorphisms $L\rightarrow M,L\rightarrow \overline{M}$ are
given by the homomorphisms
\begin{equation*}
\lambda :M\otimes \overline{M}\rightarrow M\text{ and }\lambda ^{\prime
}:M\otimes \overline{M}\rightarrow \overline{M}
\end{equation*}%
and $\partial :C\rightarrow M\cap \overline{M}.$ The $h$-map of the crossed
square is given by
\begin{equation*}
h(m,\overline{n})=i(m\otimes \overline{n})
\end{equation*}%
for $m,n\in M.$
\end{prop}

\begin{pf}
This comes by direct verification using the universal properties of tensors
and coproducts.
\end{pf}

\bigskip

\textbf{Remark}: For future applications it is again important to
note that the result is not dependent on the crossed square being
\emph{totaly} free. If $M\rightarrow R$ is any pre-crossed module,
one can form the corner
\begin{equation*}
\xymatrix{&M\ar[d]^{\nu}\\
\overline{M}\ar[r]_-{\mu}&M\rtimes R}
\end{equation*}
complete it to a crossed square via $M\otimes \overline{M}$ and then add in $%
\mathbf{S_{3}}\rightarrow M.$ Nowhere does this use freeness of
$M\rightarrow R$.

\begin{cor}
\label{c} Let \textbf{E}$^{(1)}$ be the 1-skeleton of a free
simplicial algebra. Given the free crossed square
$\mathbf{M}(\mathbf{E}^{(1)},2)$ described above, then
\begin{equation*}
NE_{2}^{(1)}/\partial _{3}NE_{3}^{(1)}\cong \mbox{\rm
Ker}d_{1}^{1}\otimes _{E_{1}}\mbox{\rm Ker}d_{0}^{1}.
\end{equation*}
\end{cor}

\begin{pf}
In the 1-skeleton of a free simplicial algebra \textbf{E}$^{(1)},$ the set
\textbf{S}$_{3}$ is empty. Thus this is clear from the previous proposition.
\end{pf}

\textbf{Remark}: If we set $M=\ker d_{0}^{1}=NE_{1}^{(1)}$, then the
identification given by the Corollary gives
\begin{equation*}
NE_{2}^{(1)}/\partial _{3}(NE_{3}^{(1)})\cong M\otimes \overline{M}.
\end{equation*}%
This uses the fact that $\ker d_{0}^{1}$ and $\ker d_{1}^{1}$ are linked via
the map sending $m$ to $\left( m-s_{0}d_{1}m\right) $ for $m\in \ker
d_{0}^{1}.$ The $h$-map
\begin{equation*}
h:M\times \overline{M}\longrightarrow NE_{2}^{(1)}/\partial _{3}NE_{3}^{(1)}
\end{equation*}%
given by
\begin{equation*}
h(x,\overline{y})=s_{1}x(s_{1}y-s_{0}y)+\partial _{3}NE_{3}^{(1)}.
\end{equation*}%
But this is also $h(x,\overline{y})=x\otimes \overline{y}$ . Thus
\begin{equation*}
x\otimes \overline{y}=s_{1}x(s_{1}y-s_{0}y)+\partial
_{3}NE_{3}^{(1)}
\end{equation*}%
under the identification via the isomorphism of the above corollary.

This explains the mysterious formula of \cite{ap2} in the discussion
before Proposition 2.6 of that paper.

\section{Applications}

\subsection{2-crossed complexes}

A notion of $2$-crossed complex of commutative algebras is defined
by Grandjean and Vale in \cite{gv}. We have considered freeness
conditions in \cite{ap3} and this generalises easily to 2-crossed
complexes.

A \emph{2-crossed complex} of commutative algebras is a sequence
of $k$-algebras
\begin{equation*}
\xymatrix{C:\cdots \ar[r]^-{\partial_{n+1}}&C_n
\ar[r]^-{\partial_n}& C_{n-1}\ar[r]^-{\partial_{n-1}}&\cdots
\ar[r]^{\partial_3}&C_2 \ar[r]^{\partial_2 }&C_1
\ar[r]^{\partial_1}&C_0}
\end{equation*}
together with a 2-crossed module structure given by the pairing
\begin{equation*}
\left\{ \quad \otimes \quad \right\} :C_{1}\otimes
_{C_{0}}C_{1}\longrightarrow C_{2}
\end{equation*}
such that

(i) $C_{n}$ is an $A$-module for $n\geq 3$ with $A=C_{0}/\partial
_{1}(C_{1});$

(ii) $C_{0}$ acts on $C_{n},$ $n\geq 1,$ the action of $\partial _{1}(C_{1})$
being trivial on $C_{n}$ for $n\geq 3;$

(iii) each $\partial _{n}$ is an $A$-module homomorphism and $\partial
_{n}\partial _{n+1}=0$ for all $n\geq 1.$

\bigskip

Note that $K=$ Ker $\partial _{2}$ is an $C_{0}/\partial _{1}C_{1}$-module
as $\partial _{2}$ is a crossed module.

The notion of a morphism for 2-crossed complexes should be clear.
Such a morphism will be a morphism `chain complexes of algebras'
restricting to a morphism of 2-crossed modules on the bottom three
terms and compatible with the action. This gives the category,
\textbf{2-CrsComp}, of 2-crossed complexes and morphisms between
them.
\begin{prop}
There is a functor%
\begin{equation*}
\mathbf{C}^{(2)}:\mathbf{SimpAlg}\longrightarrow
\mathbf{2\text{-}CrsComp}.
\end{equation*}
\rm{(We will usually omit the superfix $(2)$ writing simply
$\mathbf{C}$ for this.)}
\end{prop}

\begin{pf}
\bigskip Given a simplicial algebra $\mathbf{E},$ define
\begin{equation*}
C_{n}=
\begin{cases}
NE_n & \text{for } n=0,1;\\
\dfrac{NE_{2}}{\partial _{3}(NE_{3}\cap D_{3})} &   \text{for }  n=2; \\
\dfrac{NE_{n}}{(NE_{n}\cap D_{n})+\partial _{n+1}(NE_{n+1}\cap
D_{n+1})} & \text{for }  n\geq 3,%
\end{cases}%
\end{equation*}%
with $\partial _{n}$ induced by the differential of \textbf{NE}.
Note that the bottom three terms (for $n=0,1$ and $2$) form the
2-crossed module considered in \cite{ap3} and that for $n\geq 3$ the algebras are all $%
C_{0}/\partial _{1}(C_{1})$-modules, since in these dimensions
$C_{n}$ is the same as the corresponding crossed complex term (cf.
\cite{ap1}). The only thing remaining is to check that $\partial
_{2}\partial _{3}$ is trivial which is straightforward.
\end{pf}

Since
\begin{equation*}
\mathbf{C^{(2)}}(\mathbf{E})_{2}=NE_{2}/\partial _{3}(NE_{3}\cap
D_{3}),
\end{equation*}
the same formula as that for $\mathbf{C(E},2)_{\langle 2 \rangle}$,
we obtain the following result.

\begin{cor}
If $\mathbf{E}^{(1)}$ is the $1$-skeleton of a free simplicial
algebra $\mathbf{E}$ then the $2$-crossed complex of
$\mathbf{E}^{(1)}$ satisfies
\begin{equation*}
\mathbf{C^{(2)}}(\mathbf{E}^{(1)})_{2}\cong
\text{Ker}d_{1}^{1}\otimes \text{Ker}d_{0}^{1}.
\end{equation*}%
Moreover%
\begin{equation*}
\mathbf{C^{(2)}}(\mathbf{E}^{(2)})_{2}\cong \left( \left( \text{Ker}d_{1}^{1}\otimes \text{%
Ker}d_{0}^{1}\right) \sqcup C\right) /\sim
\end{equation*}%
where as in Proposition $5.1,$ this $C$ is a free crossed module on the
\textquotedblleft new generators\textquotedblright\ in dimension $2.$
\end{cor}

\begin{lem}
If  $\mathbf{E}$ is a simplicial resolution of $B$  then, for $k\geq
3,\mathbf{C^{(2)}}(\mathbf{E}^{(2)})_{k}$ a free $B$-module on the
given data.
\end{lem}

To sum up we have the following result.

\begin{thm}
The `step-by-step' construction of a simplicial resolution of an
algebra $B$ gives a  `step-by-step' construction of a 2-crossed
resolution of $B$ via the 2-crossed complex construction
$\mathbf{C}^{(2)}$.
\end{thm}

\subsection{`Quadratic' analogues of the cotangent complex?}

In this final section we take the $1$-skeleton of a simplicial
algebra and see how it relates to other algebraic construction,
such as Andr\'{e}-Quillen homology and a squared complex form of
the cotangent complex. Here our results are less conclusive than
we would like.

\textbf{Remark}: We will assume that rings and algebras are
Noetherian for convenience and thus that ideals are finitely
generated. This means $\mathbf{S_2, S_3}$ etc will all be finite.

From Proposition 4.1 we have the free crossed square
$$
 \left(
\begin{array}{cc}
\xymatrix{\dfrac{(R[s_0({\bf S_2})]^+[s_1({\bf S_2})])\cap (s_0({\bf
S_2})-s_1({\bf S_2}))}{P_2} \ar[d]\ar[r]&\overline{R^+[{\bf
S_2}]}\ar[d]\\
   R^+[{\bf S_2}]\ar[r]&R[{\bf S_2}]}
\end{array}
\right)
$$
so using corollary 5.2 there is the following isomorphism:
$$
\dfrac{(R[s_0({\bf S_2})]^+[s_1({\bf S_2})])\cap (s_0({\bf
S_2})-s_1({\bf S_2}))}{P_2} \cong R^{+}[\mathbf{S_2}]\otimes
\overline{R^{+}[\mathbf{S_2}]}.
$$
Thus the free crossed square becomes
$$
 \left(
\begin{array}{cc}
\xymatrix{ R^{+}[\mathbf{S_2}]\otimes
\overline{R^{+}[\mathbf{S_2}]}\ar[d]\ar[r]&\overline{R^+[{\bf
S_2}]}\ar[d]\\
   R^{+}[{\bf S_2}]\ar[r]&R[{\bf S_2}].}
\end{array}
\right)
$$
In section 4 we saw that there is a 2-crossed module
$$
\mathbb{X}:\xymatrix{&R^{+}[\mathbf{S_2}]\otimes
\overline{R^{+}[\mathbf{S_2}]}\ar[r]^{\partial_3}&R^{+}[\mathbf{S_2}]\rtimes
\overline{R^{+}[\mathbf{S_2}]}\ar[r]^-{\partial_2}&R[\mathbf{S_2}]}
$$
where
$$
\partial_3(x\otimes y)=(-\lambda(x\otimes y),\lambda'(x\otimes
y))
$$
and
$$
\partial_2(x,y)=\mu(x)+\mu'(x).
$$
The axioms of a squared complex ensure that $\partial_3$ and
$\partial_2$ are homomorphisms and $\partial_3$ is a module.

The 2-crossed complex $\mathbf{C^{(2)}(E^{(2)})}$ has a smaller
2-crossed module at its base namely
$$
\mathbb{Y}:\xymatrix{&R^{+}[\mathbf{S_2}]\otimes
\overline{R^{+}[\mathbf{S_2}]}\ar[r]&R^{+}[\mathbf{S_2}]\ar[r]^-{\partial_2}&R}
$$
and it is important to compare the two. In fact there is a split
epimorphism from $\mathbb{X}$ to $\mathbb{Y}$ with kernel
$$
\xymatrix{0\ar[r]&R^{+}[\mathbf{S_2}]\ar[r]^{=}&R^{+}[\mathbf{S_2}]}
$$
which has, of course, trivial homotopy. Thus $\mathbb{X}$ and
$\mathbb{Y}$ encode the same information about the presentation of
$R/I$, $I=$Im $\partial_2$.

Crossed complexes form a category \textbf{Crs} which can be
considered as a full subcategory of both the categories of 2-crossed
complexes and of squared complexes. In the case of 2-crossed
complexes, any crossed complex
$$
\mathbb{C}:\xymatrix{\cdots\ar[r]&C\ar[r]&M\ar[r]&R}
$$
yields a 2-crossed complex with the same terms at each level and
with trivial Peiffer lifting $\{\otimes\}:M\otimes_{R} M\rightarrow
C$ whilst considered as a squared complex we get $\mathbb{C}$ yields
$$
\xymatrix{\cdots\ar[r]&C\ar[r]\ar[d]&M\ar[d]\\
&0\ar[r]&R.}
$$
In both cases higher dimensional terms are left unchanged. Both
these inclusions have left adjoints, i.e., the embeddings give
reflective subcategories. The proofs are quite easily (and will be
given elsewhere).

The functors from \textbf{SimpAlg} to \textbf{2-CrsComp} and
\textbf{SqComp} used above, when composed with the reflections to
\textbf{Crs} yield the associated crossed complex functor mentioned
earlier (see \cite{ap1}).

Finally the category of chain complexes over $R/I$ embeds as a
reflexive subcategory of \textbf{Crs} and the reflection  sends a
crossed resolution to the (intermediate stage of the) cotangent
complex (see \cite{ap1} and \cite{porter2}). Thus given a simplicial
resolution of $R/I$, constructed as in \cite{andre} by a
step-by-step method, the 2-crossed and squared resolutions it gives
can be considered as `quadratic' analogues of the cotangent complex,
in the same way that the crossed complex is a `linear' homotopy
analogue of the `homological' cotangent complex. (Here we are using
`quadratic' and `linear' in the analogous way to that used by Baues
in \cite{b1} for the group based theory.)

Given this it is of interest to study the complex $\mathbb{X}$ (or
equivalently  $\mathbb{Y}$) and their analogues when $\mathbf{S_3}$
information is added in. Here we have no definitive results, only
problems.

The idea will be to try to provide algorithms for calculating and
thus controlling, the kernel of $\partial_3$  in  $\mathbb{X}$ (or
equivalently  $\mathbb{Y}$). We know these give
$\pi_3(\mathbf{E}^{(1)})$ and it is hoped that if these algorithms
worked, they would allow an analysis of $\pi_3(\mathbf{E}^{(2)})$
and thus to study the effect of adding in $\mathbf{S_3}$ information
to the higher terms of the simplicial resolution. As yet we are not
sure if a general analysis will be possible or whether it will be
necessary to limit ourselves to specific classes of example, using,
for instance, methods from Gr\"{o}bner base theory.

\bigskip $
\begin{array}{llllllll}
\text{Z. Arvasi} &  &  &  &  &  & \text{E. Ulualan} &  \\
\text{Eski\c{s}ehir Osmangazi University}, &  &  &  &  &  & \text{Dumlup\i
nar University} &  \\
\text{Science and Art Faculty} &  &  &  &  &  & \text{Science and Art Faculty%
} &  \\
\text{Department of Mathematics and Computer Science} &  &  &  &  &  & \text{Mathematics Department}
&  \\
\text{26480, Eski\c{s}ehir, TURKEY} &  &  &  &  &  & \text{K\"{u}tahya,
TURKEY} &  \\
\text{e-Mail: zarvasi@ogu.edu.tr, } &  &  &  &  &  & \text{%
eulualan@dumlupinar.edu.tr} &
\end{array}%
$

\end{document}